\documentclass[11pt]{article}
\usepackage{amsmath,amssymb}

\newtheorem{theo}{Theorem}[section]

\newtheorem{lemma}[theo]{Lemma}
\newtheorem{claim}[theo]{Claim}

\newtheorem{conj}[theo]{Conjecture}

\newcommand{\ignore}[1]{}

\newcommand\be{\begin{eqnarray}}
\newcommand\ee{\end{eqnarray}}

\begin{document}

\title{The domatic number of regular and almost regular graphs}

\author{Raphael Yuster\\
Department of Mathematics\\ University of Haifa at Oranim\\ Tivon 36006,
Israel.\\ e-mail: raphy@research.haifa.ac.il}

\date{}

\maketitle

\begin{abstract}
The domatic number of a graph $G$, denoted $dom(G)$, is the maximum possible
cardinality of a family of disjoint sets of vertices of $G$, each set being a
dominating set of $G$. It is well known that every graph without isolated
vertices has $dom(G) \geq 2$. For every $k$, it is known that there are
graphs with minimum degree at least $k$ and with $dom(G)=2$. In this paper we
prove that this is not the case if $G$ is $k$-regular or {\em almost}
$k$-regular (by ``almost'' we mean that the minimum degree is $k$ and the maximum degree is at most $Ck$
for some fixed real number $C \geq 1$). In this case we prove that $dom(G) \geq
(1+o_k(1))k/(2\ln k)$. We also prove that the order of magnitude $k/\ln k$ cannot be improved.
One cannot replace the constant 2 with a constant smaller than 1. The proof uses the so called
{\em semi-random method} which means 
that combinatorial objects are generated via repeated applications of
the probabilistic method; in our case iterative applications of the 
Lov\'asz Local Lemma.
\end{abstract}

\setcounter{page}{1}

\section{Introduction}
All graphs considered here are finite, undirected and simple. For standard
graph-theoretic terminology the reader is referred to \cite{Bo}. A subset $D$
of vertices in a graph $G$ is a {\it dominating set} if every vertex not in
$D$ has a neighbor in $D$. The {\em domatic number} of a graph $G$, denoted
$dom(G)$, is the maximum number of colors in a (not necessarily proper)
vertex coloring of $G$, where each color class is a dominating set. The
theory of domination and the domatic number are well studied areas in graph
theory and theoretical computer science. The two books \cite{HaHeSl1,HaHeSl2} present most of the known
results in domination theory. The domatic number was first defined in \cite{CoHe}.

Clearly, every graph has a dominating set. Thus, let $\gamma(G)$ denote the
minimum possible cardinality of a dominating set. It is an easy observation
that every graph without isolated vertices has $dom(G) \geq 2$. Simply take a
dominating set $D$ with $|D|=\gamma(G)$ and notice that $V(G) \setminus D$ is
also a dominating set. On the other hand, deciding whether $dom(G) > 2$ is an
NP-Complete problem \cite{GaJoTa}. It is well known \cite{Lo} that if a graph
has high minimum degree, then $\gamma(G)$ is small. This, however, does not
necessarily mean that $dom(G)$ is large. Zelinka \cite{Ze} proved that for
every $k$, there are graphs with minimum degree $k$ and with
$dom(G)=2$. In his examples, there are always (relatively) few vertices with
very high degree. This is unavoidable. In this paper we show that if we only
consider the class of graphs with minimum degree $k$ and maximum degree at
most $Ck$, for some constant $C \geq 1$ (in particular, regular graphs), then, in fact, $dom(G)$ is
guaranteed to be quite large.

Before we present our main result we need a definition. Let $k$ be a
positive integer and let $C \geq 1$ be a real number. A graph $G$ is called
{\em $(k,C)$-regular} if $\delta(G)=k$ and $\Delta(G) \leq kC$. In
particular, a $(k,1)$-regular graph is a $k$-regular graph. Let $f(k,C)$
denote the minimum possible value of $dom(G)$ taken over all $(k,C)$-regular
graphs. For example, $f(2,C)=2$ as seen
be any cycle whose number of vertices is not divisible by 3. Also, $ f(3,C)=2$ as seen by
the 3-regular graph with 8 vertices consisting of a Hamiltonian cycle and four edges connecting
antipodal vertices of the cycle. In general, determining $f(k,C)$ precisely seems to be
a very difficult task. Our main result is summarized in the following two theorems:
\begin{theo}
\label{t11}
Let $C \geq 1$ be a fixed real number. Then,
$$
f(k,C) \geq \frac{k}{2\ln k}(1+o_k(1)).
$$
\end{theo}
\begin{theo}
\label{t12}
Let $C > 1$ be a fixed real number. Then,
$$
f(k,C) \leq \frac{k}{\ln k}(1+o_k(1)).
$$
\end{theo}
Although Theorem \ref{t11} and Theorem
\ref{t12} show, in particular, that the order of magnitude of $f(k,C)$ is
$k/\ln k$ for any fixed $C > 1$, the constants are worthy of
investigation. In fact, although the proof of Theorem \ref{t11} is
significantly more difficult than the proof of Theorem \ref{t12}
we conjecture that the latter is the correct answer
\begin{conj}
\label{c12}
Let $C \geq 1$ be a fixed real number. Then,
$$
f(k,C) = \frac{k}{\ln k}(1+o_k(1)).
$$
\end{conj}

In the next two sections we present the proofs of Theorem \ref{t11} and Theorem
\ref{t12}. The proof of theorem \ref{t11} demonstrates the so called semi-random
method (an unofficial term the author heard several times by researchers in the area
of probabilistic methods in combinatorics). In this method,
combinatorial objects (such as graph colorings) are generated via repeated applications of
the probabilistic method; in the case of Theorem \ref{t11}, iterative applications of the 
Lov\'asz Local Lemma. The proof of Theorem \ref{t12} is demonstrated by exhibiting an appropriate random graph.

A final note: Proving an analog of Theorem \ref{t11} with the constant 3 instead of 2 is a significantly easier task.
In fact, a naive application of the Local Lemma does the job. Assume you have about $t=\lfloor k/(3 \ln k) \rfloor$ colors.
Let each vertex choose a random color independently with uniform distribution.
Let $A_{v,i}$ denote the event that the vertex $v$ misses the color $i$ in its closed neighborhood.
The probability of $A_{v,i}$ is less than $1/k^3$. But $C^2k^3/\ln k$ is also an upper bound for the number of
events which $A_{v,i}$ depends upon, since if $u$ and $v$ are at distance 3 or more from each other, $A_{v,i}$ is independent
of $A_{u,j}$ (they do not have common neighbors). The total number of vertices at distance 2 is
at most $kC+k^2C^2$ and there are $t$ colors, so the dependency digraph has maximum degree
$(kC+k^2C^2)k/(3 \ln k)$. Now the conditions of the Local Lemma hold, and thus with positive probability
no $A_{v,i}$ holds. Hence each color class is a dominating set. This naive approach fails, of course, for any
constant smaller than 3, so additional ideas must be sought.

\section{Proof of Theorem \ref{t11}}
In the proof of Theorem \ref{t11} we need to use the
Lov\'asz Local Lemma \cite{ErLo}. Here it is,
following the notations in \cite{AlSp} (which also contains a simple proof of
the lemma). Let $A_1,\ldots,A_n$ be events in an arbitrary probability space.
A directed graph $D=(V,E)$ on the set of vertices $V=[n]$ is called a {\em
dependency digraph} for the events $A_1,\ldots,A_n$ if for each $i$,
$i=1,\ldots,n$, the event $A_i$ is mutually independent of all the events
$\{A_j ~ : ~ (i,j) \notin E\}$.
\begin{lemma}[The Local Lemma, symmetric version]
\label{l21}
Let $A_1,\ldots,A_n$ be events in an arbitrary probability space and let
$D=(V,E)$ be a corresponding dependency digraph.
If the maximum outdegree in $D$ is at most $d \geq 1$ and each $A_i$ has $\Pr[A_i] \leq
p$ and $p(d+1) \leq 1/e$ then with positive probability no event $A_i$ holds.
\end{lemma}

\noindent
{\bf Proof of Theorem \ref{t11}:}\,
Let $C \geq 1$ be fixed, and let $\epsilon > 0$.
in order to avoid cluttered computations we shall assume,
wherever necessary, that $k$ is sufficiently large as a function
of $C$ and $\epsilon$ only.
Let $k$ be sufficiently large such that there is an integer between
$k/((2+\epsilon)\ln k)$ and $k/((2+\epsilon/2)\ln k)$.
Thus, for some $\epsilon/2 \leq \gamma \leq \epsilon$,
The number $t=k/((2+\gamma)\ln k)$ is an integer.
Let $G=(V,E)$ be a $(k,C)$-regular graph.
We need to show that $dom(G) \geq t$. This will show, in
particular, that $dom(G)\geq k/((2+\epsilon)\ln k)$ and consequently, $f(k,C)
\geq (k/2\ln k)(1+o_k(1))$.

Assume that we have the set of colors $\{0, 1,\ldots,t\}$.
We call color $0$ the {\em transparent} color.
In the first phase of the proof we color the vertices using all colors such that
certain very specific properties hold. In the second phase we recolor the
vertices that received the transparent color in the first phase using only
the non-transparent colors and show that we can do it carefully enough such that each non-transparent color class
(after the second phase) is a dominating set.

We begin with a description of the first phase.
Our goal in the first phase is to achieve a coloring with the following properties:
\begin{lemma}
\label{l22}
There exists a coloring of $G$ with the colors $\{0,1,\ldots,t\}$ such that the following conditions hold:
\begin{enumerate}
\item
Every vertex has at least $k\gamma/(4(2+\gamma))$ neighbors with transparent color.
\item
Every vertex has at most 4 non-transparent colors missing from its (open) neighborhood.
\item
Put $z=\lceil 12/\gamma \rceil$. For each $v \in V$, and for each sequence of
$z$ {\bf distinct} non-transparent colors $c_1,\ldots,c_z$ and for each
sequence of $z$ {\bf distinct} neighbors of $v$ denoted $u_1,\ldots,u_z$, at least one $u_i$
has a neighbor colored $c_i$.
\end{enumerate}
\end{lemma}
{\bf Proof:}\,
We let each vertex $v \in V$ choose one
color from $\{0, 1,\ldots,t\}$ randomly.
The probability to choose color $i$ is $p=(2+\gamma/2)\ln k/k$ for $i=1,\ldots,t$
and the probability to choose the transparent color is, therefore, $q=1-pt=\gamma/(2(2+\gamma))$.
Let $A_v$ denote the event that $v$ has less than $kq/2$ neighbors colored with the transparent color.
Let $B_v$ denote the event that $v$ has more than 4 non-transparent colors missing from its
neighborhood. Let $C_v$ denote the event that $v$ has $z$ neighbors $u_1,\ldots,u_z$
and there exist $z$ distinct non-transparent colors $c_1,\ldots,c_z$, such that $c_i$ is
missing from the neighborhood of $u_i$ for each $i=1,\ldots,z$.
Thus, we need to show that with positive probability, none of the
$3|V|$ events $A_v$,$B_v$ and $C_v$, for each $v \in V$, hold.
The following three claims provide upper bounds for the probabilities of the events $A_v$,$B_v$ and $C_v$,
respectively.
\begin{claim}
\label{c1}
$\Pr[A_v] < 1/k^5$.
\end{claim}
{\bf Proof:}\,
Let $X_v$ denote the random variable counting the number of transparent neighbors of $v$.
The expectation of $X_v$ is
$E[X_v]=d_vq \geq kq$, where $d_v$ denotes the degree of $v$.
Since each vertex chooses its color independently we have by the
most common Chernoff inequality (cf. \cite{AlSp})
$$
\Pr\left[A_v\right] =\Pr\left[X_v < \frac{kq}{2}\right] \leq \Pr\left[X_v < \frac{E[X_v]}{2}\right] <
$$
$$
e^{-2(E[X_v]/2)^2/d_v}=e^{-d_v^2q^2/(2d_v)}=e^{-d_vq^2/2}<e^{-kq^2/2} << \frac{1}{k^5}.
$$
(In the final inequality we used the fact that $q$ is a constant depending on $\gamma$ and that $k$
is sufficiently large).
\begin{claim}
\label{c2}
$\Pr[B_v] < 1/k^5$.
\end{claim}
{\bf Proof:}\,
Fix 5 distinct non-transparent colors. The probability that none of them appear in the
neighborhood of $v$ is precisely $(1-5p)^{d_v}$.
Now,
$$
(1-5p)^{d_v} \leq (1-5p)^k = \left(1-\frac{5(2+\frac{\gamma}{2})\ln k}{k}\right)^k < \frac{1}{k^{10+2.5\gamma}}.
$$
As there are ${t \choose 5}$ possible sets of 5 distinct non-transparent colors
we get that
$$
\Pr[B_v] < {t \choose 5}\frac{1}{k^{10+2.5\gamma}} < \frac{1}{k^{5+2.5\gamma}} < \frac{1}{k^5}.
$$
\begin{claim}
\label{c3}
$\Pr[C_v] < 1/k^5$.
\end{claim}
{\bf Proof:}\,
For a vertex $u$ and a color $c$ let $n(u,c)$ denote the number of neighbors of
$u$ colored $c$.
Fix a set of $z$ distinct non-transparent colors $\{c_1,\ldots,c_z\}$ and $z$ distinct neighbors
of $v$, $\{u_1,\ldots,u_z\}$. We begin by computing the probability that
for each $i=1,\ldots,z$, $c_i$ does not appear in the neighborhood of $u_i$
(i.e. $n(u_i,c_i)=0$).
Denoting this probability by $\rho=\rho(v,u_1,\ldots,u_z,c_1,\ldots,c_z)$ we clearly have:
$$
\rho=\Pr[n(u_1,c_1)=0] \cdot \Pr[n(u_2,c_2)=0 ~ | ~ n(u_1,c_1)=0] \cdot \ldots
$$
$$
\ldots \cdot \Pr[n(u_z,c_z)=0 ~ | ~ n(u_1,c_1)=0 ~\land ~\cdots ~\land ~n(u_{z-1},c_{z-1})=0].
$$
For the first term we have $\Pr[n(u_1,c_1)=0]=(1-p)^{d_{u_1}}$.
For the other terms we claim that
$$
\Pr[n(u_i,c_i)=0 ~ | ~ n(u_1,c_1)=0 ~\land ~\cdots ~\land ~n(u_{i-1},c_{i-1})=0] \leq (1-p)^{d_{u_i}}.
$$
This is obvious since the knowledge that a color from $\{c_1,\ldots,c_{i-1}\}$ does not appear
in a neighbor common to $u_i$ and some $u_{i'}$ for $i' < i$ only increases the probability that $c_i$ is
in the neighborhood of $u_i$,
and hence decreases the probability that $n(u_i,c_i)=0$. To be precise, if $w_{i,1}, \ldots,
w_{i,d_{u_i}}$ are the neighbors of $u_i$, let $s(j)$ denote the size of the intersection
of $N(w_{i,j})$ with $\{u_1,\ldots,u_{i-1}\}$. Clearly $0 \leq s(j) \leq i-1 < z < t$.
The probability that $w_{i,j}$ is colored with
$c_i$ given that $n(u_l,c_l)=0$ for $l=1,\ldots, i-1$ is precisely $p/(q+(t-s(j))p)$.
Recalling that $q+tp=1$ we have $p/(q+(t-s(j))p) \geq p$. Thus,
$$
\Pr[n(u_i,c_i)=0 ~ | ~ n(u_1,c_1)=0 ~\land ~\cdots ~\land ~n(u_{i-1},c_{i-1})=0] =
\Pi_{j=1}^{d_{u_i}}\left(1-\frac{p}{q+(t-s(j))p}\right) \leq (1-p)^{d_{u_i}}.
$$
We therefore have:
$$
\rho \leq \Pi_{i=1}^{z}(1-p)^{d_{u_i}} \leq (1-p)^{kz}.
$$
There are less than $(Ck)^z$ ordered sets of $z$ distinct neighbors of $v$.
There are less than $t^z$ ordered sets of $z$ distinct non-transparent colors.
Thus,
$$
\Pr[C_v] < t^z(Ck)^z(1-p)^{kz} \leq C^zk^{2z}(1-\frac{2+(\gamma/2)\ln k}{k})^{kz} <
$$
$$
C^zk^{2z}\frac{1}{k^{z(2+\gamma/2)}} =\frac{C^z}{k^{z\gamma/2}}<\frac{C^{13/\gamma}}{k^6} < \frac{1}{k^5}.
$$ 

\noindent
Having proved $\Pr[A_v] < 1/k^5$, $\Pr[B_v] < 1/k^5$ and $\Pr[C_v] < 1/k^5$
we claim that we can use the Local Lemma to show that with positive probability none
of these events hold. Indeed, fix a vertex $v$ and let $U_v$ be the set of all vertices
at distance 5 or greater from $v$. Notice that if $u \in U_v$, the neighbors of
$v$ and their neighborhoods do not intersect the neighbors of $u$ and
their neighborhoods.
Since $A_v$ only depends on $v$ and its neighbors, $B_v$ depends only
on $v$ and its neighbors and $C_v$ only depends on $v$, its neighbors and 
the neighbors of its neighbors, we have that $A_v$ is mutually independent of
all the $3|U_v|$ events $\{A_u, B_u, C_u ~ | ~ u \in U_v\}$. Similarly $B_v$ and $C_v$ are mutually
independent of all the event $\{A_u, B_u, C_u ~ | ~ u \in U_v\}$.
Since there are at most $1+kC+kC(kC-1)+kC(kC-1)^2+kC(kC-1)^3 < k^4C^4$
vertices at distance at most 4 from $v$ (including $v$), we have that the maximum outdegree
in the dependency digraph of the $3|V|$ events is at most $3k^4C^4$.
Since $(1/k^{5}) \cdot (3k^4C^4+1) < 1/e$ we get by Lemma \ref{l21} that with positive
probability none of the events hold.
We therefore proved Lemma \ref{l22}. 

We now describe the second phase. We fix a coloring satisfying the three conditions
in the statement of Lemma \ref{l22}.
For a vertex $v$, let $F(v)$ denote the set of missing non-transparent colors from its
neighborhood. By Lemma \ref{l22} we know that $|F(v)| \leq 4$.
Now, let $S(v) = \cup_{u \in N(v)} F(u)$.
We claim that $|S(v)| \leq 4(z-1)$. To see this, notice that if $|S(v)| > 4(z-1)$ this means
that there are at least $z$ distinct neighbors of $u$, each missing a distinct color
from their neighborhood, contradicting the third condition in Lemma \ref{l22}.
In the second phase we only color the vertices that received transparent colors
in the first phase. Let $v$ be a vertex colored with the transparent color.
We let $v$ choose a random color from $S(v)$ with uniform distribution.
The choices made by distinct vertices are independent (In case $S(v)=\emptyset$ 
we can assign an arbitrary non-transparent color to $v$).
Let $v \in V$ be any vertex, and let $c \in F(v)$. Let $A_{v,c}$ denote
the event that after the second phase, $c$ still does not appear as a color in
a neighbor of $v$. Our goal is to show that with positive probability, none of
the events $A_{v,c}$ for $v \in V$ and $c \in F(v)$ hold. This will complete the
proof of Theorem \ref{t11}.

Let $T_v$ be the subset of neighbors of $v$ given transparent color in the
first phase. By Lemma \ref{l22} we have $|T_v| \geq k\gamma/(4(2+\gamma))$.
Assuming $c$ does not appear in the neighborhood of $v$ we have that
for each $u \in T_v$, the color $c$ appears in $S(u)$.
Hence,
$$
\Pr[A_{v,c}] = \Pi_{u \in T_v}\left(1-\frac{1}{|S(u)|}\right) < \Pi_{u \in T_v}\left(1-\frac{1}{4z}\right)<
$$
$$
\left(1-\frac{1}{4z}\right)^{k\gamma/(4(2+\gamma))}<
\left(1-\frac{\gamma}{52}\right)^{k\gamma/(4(2+\gamma))} << \frac{1}{k^3}.
$$
Now, for $v \in V$, let $U_v$ denote the set of all vertices at distance at least 3 from $v$.
Since the event $A_{v,c}$ only depends on $v$ and its (transparent) neighbors, we have
that $A_{v,c}$ is mutually independent of all the events $A_{u,c'}$ for $u \in U_v$
and $c' \in F(u)$. Since the number of neighbors at distance at most 2 from $v$ is at
most $1+kC+kC(kC-1)$, including $v$, and since $|F(u)| \leq 4$ for all $u \in V$ we
have that the outdegree in the dependency digraph of the events is at most $4(1+kC+kC(kC-1))< 5k^2C^2$.
Since $(1/k^{3}) \cdot (5k^2C^2+1) < 1/e$ we get by Lemma \ref{l21} that with positive
probability none of the events of the form $A_{v,c}$ hold.
Hence, there is a coloring with the colors $\{1,\ldots,t\}$ such that each color class
is a dominating set.

\section{Proof of Theorem \ref{t12}}
We shall take the opportunity to prove something slightly stronger than the statement of Theorem \ref{t12}.
The random graphs we shall construct to demonstrate the proof of Theorem \ref{t12}
can also have arbitrary large girth.

Trivially, $dom(G) \leq |V(G)|/\gamma(G)$.
Thus, it suffices to prove the following:
\begin{theo}
\label{t31}
Let $C > 1$ be a fixed real number, and let $g \geq 2$.
For every $\epsilon > 0$, there exists a $k_0=k_0(C,g,\epsilon)$ such that for all $k > k_0$,
there exists a $(k,C)$-regular graph $G$ with $n$ vertices and
with $girth(G) > g$ having $\gamma(G) \geq (1-\epsilon)n\ln k/k$.
\end{theo}
A weaker theorem, in which we only require the graph $G$ to have minimum degree $k$
(i.e. $C=\infty$)
and we do not care about the girth (i.e. $g=2$) follows immediately from a result
of Alon \cite{Al}.
{\bf Proof of Theorem \ref{t31}:}\,
Put $r=2g+1$. Trivially, we may assume
\begin{equation}
\label{e31}
1+\epsilon/4 < C ~ , ~ \epsilon < 1.
\end{equation}
Let $k_0$ be the minimal integer satisfying
\begin{enumerate}
\item
\begin{equation}
\label{e32}
k_0 \geq \frac{48}{\epsilon^2}.
\end{equation}
\item
\begin{equation}
\label{e33}
k_0 \geq 6g \cdot 8^{g^2}.
\end{equation}
\item
For every $k > k_0$

\begin{equation}
\label{e34}
2k^r \exp\left(-k \frac{\epsilon^2} {1024(1+\epsilon/8)}\right) < \frac{1}{3}.
\end{equation}
\item
For every $k > k_0$
\begin{equation}
\label{e35}
k^{\epsilon/4}(1-\epsilon)-r(\ln k)^2 > 1.
\end{equation}
\end{enumerate}
Let $k > k_0$ and let $n=k^r$.
Consider the random graph $G(n,p)$ where $p=(1+\epsilon/8)/k^{r-1}$.
That is, every edge appears, independently, with probability $p$.
We shall prove the following three lemmas, which, together, supply the required
result.
\begin{lemma}
\label{l31}
With probability greater than $2/3$, $G$ has minimum degree at least $k+1$ and
maximum degree at most $kC$.
\end{lemma}
\begin{lemma}

\label{l32}
With probability greater than $2/3$, $\gamma(G) \geq (1-\epsilon)n\ln k/k$.
\end{lemma}
\begin{lemma}
\label{l33}
With probability greater than $2/3$, any two cycles $C$ and $C'$ of
$G$ having  at most $g$ vertices each, are vertex disjoint.
\end{lemma}
By Lemmas \ref{l31}, \ref{l32} and \ref{l33} we know that with positive
probability there exists an $n$-vertex graph $G$ that has minimum degree at least
$k+1$, maximum degree at most $kC$,
has $\gamma(G) \geq (1-\epsilon)n\ln k/k$ and any two cycles in $G$ whose lengths
are at most $g$ are vertex-disjoint. Thus, we can delete a single edge from each
cycle whose length is at most $g$. The resulting graph is $(k,C)-regular$
(no vertex lost more than one edge), has girth greater than $g$,
and $\gamma(G)$ cannot decrease when we delete edges.
This proves Theorem \ref{t31} and, consequently, Theorem \ref{t12}. 

\noindent
{\bf Proof of Lemma \ref{l31}:}\,
The proof of this lemma is almost trivial, and is based on standard large
deviation approximations. For $v \in V(G)$, let $d_v$ denote its degree in
$G(n,p)$. $d_v$ is a random variable with the binomial distribution $B(n-1,p)$.
Thus, $E[d_v] = (n-1)p = (k^r-1)(1+\epsilon/8)/k^{r-1}$. Clearly, by (\ref{e32})
$k(1+\epsilon/8) > E[d_v] > k(1+\epsilon/12)$.
We shall use the large deviation inequality of Chernoff (cf. \cite{AlSp}) that states that
for all $a > 0$
$$
\Pr[|d_v-E[d_v]| > a] < 2\exp\left(-\frac{a^2}{2p(n-1)} + \frac{a^3}{2p^2(n-1)^2}\right).
$$
Using $a=k\epsilon/16$ and the last inequality we get, together with (\ref{e34}), that
$$
\Pr[|d_v - E[d_v]| > k\epsilon/16] < 2\exp\left( -\frac{k^2\epsilon^2}{256} \cdot \frac{1}{2k(1+\epsilon/8)}
+\frac{k^3\epsilon^3}{4096} \cdot \frac{1}{2k^2(1+\epsilon/12)^2}\right) =
$$
$$
2\exp\left(-\frac{k\epsilon^2}{512(1+\epsilon/8)}+\frac{k\epsilon^3}{8192(1+\epsilon/12)^2}\right)
< 2\exp\left(-\frac{k\epsilon^2}{1024(1+\epsilon/8)}\right) < \frac{1}{3k^r}.
$$
As there are $n=k^r$ vertices in $G$ we get that with probability greater than
$1-k^r/(3k^r) = 2/3$, all $v \in V$ satisfy $|d_v-E[d_v]| \leq k\epsilon/16$.
In particular, $d_v \geq E[d_v]-k\epsilon/16 \geq k(1+\epsilon/12-\epsilon/16) \geq k+1$,
and $d_v \leq E[d_v] + k\epsilon/16 \leq k(1+\epsilon/8+\epsilon/16) < kC$.
Thus, with probability at least 2/3, the minimum degree is at least $k+1$ and the
maximum degree is at most $Ck$. 

\noindent
{\bf Proof of Lemma \ref{l32}:}\,
Let $\epsilon/2 < \alpha < \epsilon$ be such that $t=(1-\alpha)k^{r-1}\ln k$ is an integer.
By (\ref{e32}) $\alpha$ exists.
We must show that with probability greater than $2/3$, every subset of $t$ vertices
is not a dominating set. Fix $X \subset V(G)$ with $|X|=t$. For $v \in V(G) \setminus X$,
the probability that $v$ is not adjacent to any vertex of $X$ is precisely
$(1-p)^t$. Thus, $v$ is dominated by $X$ with probability $1-(1-p)^t$.
Since the edges of $G$ are chosen independently, the probability that
$|X|$ is a dominating set is precisely $(1-(1-p)^t)^{n-t}$. As there are ${n \choose t}$
choices for $X$ it suffices to show  that ${n \choose t}(1-(1-p)^t)^{n-t} < 1/3$.
First, notice that by (\ref{e32}) and $\epsilon/2 < \alpha < \epsilon$ we have that
$$
(1-p)^t=\left(1-\frac{(1+\epsilon/8)}{k^{r-1}}\right)^{(1-\alpha)k^{r-1}\ln k} =
\left(1-\frac{(1+\epsilon/8)}{k^{r-1}}\right)^
{
\left(\frac{k^{r-1}}{1+\epsilon/8}-1\right)\frac{(1-\alpha)k^{r-1}\ln k}{\left(\frac{k^{r-1}}{1+\epsilon/8}-1\right)}
} > 
$$
$$
\left(1-\frac{(1+\epsilon/8)}{k^{r-1}}\right)^{
\left(\frac{k^{r-1}}{1+\epsilon/8}-1\right)\frac{(1-\alpha)k^{r-1}\ln k}{k^{r-1}(1-\epsilon/4)}
} >  \exp\left(-\frac{(1-\alpha)\ln k}{1-\epsilon/4}\right) > \exp\left(-\frac{(1-\epsilon/2)\ln k}{1-\epsilon/4}\right) >
$$
$$
\exp\left(-(1-\epsilon/4)\ln k\right) =\frac{1}{k^{1-\epsilon/4}}.
$$
Thus, using the fact that $n-t -k^r(1-\epsilon)=k^{r-1}(\epsilon k - (1-\alpha)\ln k) > 0$ that
follows from (\ref{e32}), and using (\ref{e35}) and the last inequality we have that
$$
{n \choose t}\left(1-(1-p)^t\right)^{n-t}  < {n \choose t}\left(1-\frac{1}{k^{1-\epsilon/4}}\right)^{n-t} <
{n \choose t}\left(1-\frac{1}{k^{1-\epsilon/4}}\right)^{k^r(1-\epsilon)} <
$$
$$
{n \choose t}\exp\left(-k^{r-1+\epsilon/4}(1-\epsilon)\right) <
(k^r)^{k^{r-1}\ln k}\exp\left(-k^{r-1+\epsilon/4}(1-\epsilon)\right) =
$$
$$
\exp\left(rk^{r-1}(\ln k)^2\right) \exp\left(-k^{r-1+\epsilon/4}(1-\epsilon)\right)=
\exp\left(-k^{r-1}(k^{\epsilon/4}(1-\epsilon)-r(\ln k)^2)\right) < \frac{1}{3}.
$$ 

\noindent
{\bf Proof of Lemma \ref{l33}:}\,
Let $F_g$ be the family of all graphs with at least $4$ vertices and at most $2g-1$
vertices, and which have more edges than vertices.
Trivially, if a graph $G$ has no element of $F_g$ as a subgraph, then all
its cycles with lengths $g$ or less are vertex disjoint.
Thus, if we can prove that
the probability that $G(n,p)$ has an element of $F_g$ as a subgraph is less than $1/3$,
we are done. 
First, notice that $|F_g| < 2g \cdot 2^{2g^2}$ as there are at most $2^{j \choose 2}$
distinct labeled graphs on $j$ vertices. Fix $H \in F_g$ and let $h$ denote the number of vertices of $H$ and
$m$ denote the number of edges of $H$. Hence, $4 \leq h \leq 2g-1$ and $m \geq h+1$.
The complete graph on $n$ vertices has less than $n^h$ labeled copies of
$H$. For each labeled copy, the probability that it belongs to $G(n,p)$ is
precisely $p^m(1-p)^{{h \choose 2} - m} \leq p^m$. Thus,
if we denote by $n(H)$ the random variable corresponding to the number of copies of $H$ in $G$
we have that the expected value of $n(H)$ is less than $p^mn^h$. Hence,
$$
E[n(H)] < p^mn^h = \frac{(1+\epsilon/8)^m}{k^{m(r-1)}}k^{rh} < (1+\epsilon/8)^{2g^2}\frac{k^{rh}}{k^{(h+1)(r-1)}}
\leq \frac{(1+\epsilon/8)^{2g^2}}{k}.
$$
In the last inequality we used the fact that $(h+1)(r-1)-rh=r-h-1=2g-h \geq 1$.
Let $n(F_g)$ denote the number of subgraphs of $G$ isomorphic to an element of $F_g$.
By linearity of expectation we get, together with (\ref{e33}), that%
$$
E[n(F_g)] \leq |F_g|\frac{(1+\epsilon/8)^{2g^2}}{k} \leq 2g \cdot 2^{2g^2}\frac{(1+\epsilon/8)^{2g^2}}{k} < \frac{1}{3}.
$$
By Markov's inequality, with probability greater than $2/3$ we have $n(F_g) = 0$.


\begin{thebibliography}{99}
\bibitem{Al} N. Alon,
{\em Transversal numbers of uniform hypergraphs}, Graphs and Combinatorics 6
(1990), 1--4.

\bibitem{AlSp} N. Alon and J.H. Spencer,
{\em The Probabilistic Method}, John Wiley and Sons Inc., New York, 1991.

\bibitem{Bo} B. Bollob\'as,
{\em Extremal Graph Theory}, Academic Press, London, 1978.

\bibitem{CoHe} E.J. Cockayne and S.T. Hedetniemi,
{\em Optimal domination in graphs}, IEEE Trans. Circuits and Systems 22
(1975), 855-857.

\bibitem{ErLo} P. Erd\"os and L. Lov\'asz,
{\em Problems and results on 3-chromatic hypergraphs and some related
questions}, Infinite and Finite Sets (A. Hajnal et al., eds.), North-Holland,
Amsterdam (1975), 609-628.

\bibitem{GaJoTa} M.R. Garey, D.S. Johnson and R.E. Tarjan,
unpublished results, 1976.

\bibitem{HaHeSl1} T. Haynes, S.T. Hedetniemi and P. Slater,
{\em Domination in Graphs: The Theory}, Marcel Dekker Publishers, New York,
1997.

\bibitem{HaHeSl2} T. Haynes, S.T. Hedetniemi and P. Slater,
{\em Domination in Graphs: Selected Topics}, Marcel Dekker Publishers, New
York, 1997.

\bibitem{Lo} L. Lov\'asz,
{\em On the ratio of optimal and integral fractional covers}, Disc. Math. 13
(1975), 383--390.

\bibitem{Ze} B. Zelinka,
{\em Domatic number and degrees of vertices of a graph}, Math. Slovaca 33
(1983), 145--147.

\end{thebibliography}
\end{document}